\documentclass[11pt, DIV10,a4paper]{article}

\usepackage{natbib}
\newcommand {\ctn}{\citet} 
\usepackage{graphicx,subfigure,amsmath,latexsym,amssymb}
\usepackage{float,epsfig,multirow,rotating,times}
\usepackage{upgreek,wrapfig}
\usepackage{comment}

\newcommand{\bphi}{\boldsymbol{\phi}}

\newcommand{\bSigma}{\boldsymbol{\Sigma}}

\newcommand{\bmu}{\boldsymbol{\mu}}

\newcommand{\bV}{\boldsymbol{V}}

\newcommand{\bU}{\boldsymbol{U}}

\newcommand{\bone}{\boldsymbol{1}}

\newtheorem{theorem}{Theorem}

\numberwithin{equation}{section}
\numberwithin{algo}{section}
\numberwithin{table}{section}
\numberwithin{figure}{section}

\usepackage{setspace}
\usepackage[mathscr]{euscript}
\usepackage[margin=1in]{geometry}
\begin{document}

\normalsize

\title{\vspace{-0.8in}
On Bayesian Asymptotics in Stochastic Differential Equations with Random Effects}
\author{Trisha Maitra and Sourabh Bhattacharya\thanks{
Trisha Maitra is a PhD student and Sourabh Bhattacharya 
is an Assistant Professor in
Interdisciplinary Statistical Research Unit, Indian Statistical
Institute, 203, B. T. Road, Kolkata 700108.
Corresponding e-mail: sourabh@isical.ac.in.}}
\date{\vspace{-0.5in}}
\maketitle%

\begin{abstract}
\ctn{Maud12} investigated asymptotic properties of the maximum likelihood estimator
of the population parameters of the random effects associated with $n$ independent 
stochastic differential equations ($SDE$'s) assuming that the $SDE$'s are independent and identical ($iid$).

In this article, we consider the Bayesian approach to learning about the population parameters, and prove
consistency and asymptotic normality of the corresponding posterior distribution in the $iid$ set-up
as well as when the $SDE$'s are independent but non-identical.
%
\\[2mm]
 
\end{abstract}

\section{Introduction}
\label{sec:intro}

Mixed effects models are appropriate when dealing with data sets consisting of variability between subjects
and also within subjects, with respect to time. Although a great deal of work on mixed effects models 
exists in the statistical literature, mixed effects models where within subject variability is modeled via
stochastic differential equations ($SDE$'s) are relatively rare. For a relatively short but comprehensive
review we refer the reader to \ctn{Maud12}, who also undertake theoretical and asymptotic investigation of a class of
$SDE$-based mixed effects models
having the following form:
for $i=1,\ldots,n$,
\begin{equation}
d X_i(t)=b(X_i(t),\phi_i)dt+\sigma(X_i(t))dW_i(t),
\label{eq:sde1}
\end{equation}
where, for $i=1,\ldots,n$, $X_i(0)=x^i$ is the initial value of the stochastic process $X_i(t)$, which 
is assumed to be continuously observed on the time interval $[0,T_i]$; $T_i>0$ assumed to be known.
The function $b(x,\varphi)$ is a known, real-valued function on $\mathbb R\times\mathbb R^d$ 
($\mathbb R$ is the real line and $d$ is the dimension); this function is known as the drift function.
The function $\sigma:\mathbb R\mapsto\mathbb R$ is the known diffusion coefficient.
In the context of statistical modelling, $X_i(\cdot)$ models the $i$-th individual.
The $SDE$'s given by (\ref{eq:sde1}) are driven by independent standard Wiener processes $\{W_i(\cdot);~i=1,\ldots,n\}$, 
and $\{\phi_i;~i=1,\ldots,n\}$, which are to be interpreted as the random effect parameters associated
with the $n$ individuals, are assumed to be independent of the Brownian motions and
independently and identically distributed ($iid$) random variables with common distribution
$g(\varphi,\theta)d\nu(\varphi)$. Here $g(\varphi,\theta)$ 
is a density with respect to a dominating measure
on $\mathbb R^d$, for all $\theta$, where 
$\theta\in\Omega\subset\mathbb R^p$ ($p\geq 2d$) is the unknown parameter of interest, which is to be estimated. 
\ctn{Maud12} impose regularity conditions that ensure existence of solutions of (\ref{eq:sde1}). The conditions, which
are also adopted by us, are as follows.
\begin{itemize}
\item[(H1)]
\begin{enumerate}
\item[(i)] The function $(x,\varphi)\mapsto b(x,\varphi)$ is $C^1$ (differentiable with continuous first derivative)
on $\mathbb R\times\mathbb R^d$, and such that there exists $K>0$ so that $$b^2(x,\varphi)\leq K(1+x^2+|\varphi|^2),$$
for all $(x,\varphi)\in\mathbb R\times\mathbb R^d$.
\item[(ii)] The function $\sigma(\cdot)$ is $C^1$ on $\mathbb R$ and
$$\sigma^2(x)\leq K(1+x^2),$$
for all $x\in\mathbb R$.
\end{enumerate}
\item[(H2)] Let $X^{\varphi}_i$ be associated with the SDE of the form (\ref{eq:sde1}) with drift function $b(x,\varphi)$.
Also letting $Q^{x^i,T_i}_{\varphi}$ denote the joint distribution of $\left\{X^{\varphi}_i(t);~t\in [0,T_i]\right\}$, it is
assumed that for $i=1,\ldots,n$, and for all $\varphi,\varphi'$, the following holds:
\[ 
Q^{x^i,T_i}_{\varphi}\left(\int_0^{T_i}\frac{b^2\left(X^{\varphi}_i(t),\varphi'\right)}{\sigma^2(X^{\varphi}_i(t))}dt<\infty\right)
=1.
\]
\item[(H3)] For $f=\frac{\partial b}{\partial\varphi_j},~j=1,\ldots,d$, there exist $c>0$ and some $\gamma\geq 0$ such that
\[
\underset{\varphi\in\mathbb R^d}{\sup}\frac{\left|f(x,\varphi)\right|}{\sigma^2(x)}\leq c\left(1+|x|^{\gamma}\right).
\]
\end{itemize}

In this article, we consider $d=1$, that is, we assume one-dimensional random effects, so that $\varphi\in\mathbb R$.
Moreover, as in \ctn{Maud12}, for statistical inference we assume that $b(x,\phi_i)$ is linear in $\phi_i$; in other
words, $b(x,\phi_i)=\phi_ib(x)$. 
Under this assumption, (H3) is not required; see \ctn{Maud12} and \ctn{Maitra14a}. 
Following \ctn{Maitra14a} we further assume that
\begin{itemize}
\item[(H1$^\prime$)] 
$b(\cdot)$ and
$\sigma(x)$ are $C^1$ on $\mathbb R$ 
satisfying $b^2(x)\leq K(1+x^2)$  
and $\sigma^2(x)\leq K(1+x^2)$
for all $x\in\mathbb R$, for some $K>0$. 
\item[(H2$^\prime$)] Almost surely for each $i\geq 1$, 
\[
\int_0^{T_i}\frac{b^2(X_i(s))}{\sigma^2(X_i(s))}ds<\infty.
\]
\end{itemize}
As in \ctn{Maud12} and \ctn{Maitra14a} here 
we assume that $\phi_i$ are normally distributed implying for $k\geq 1$, $E|\phi_i|^{2k}<\infty$ so that the following holds for all $T>0$ (see \ctn{Maud12}):
\begin{equation}
\underset{t\in[0,T]}{\sup}~E\left[X_i(t)\right]^{2k}<\infty.
\label{eq:finite_sup}
\end{equation}



In fact, the linearity assumption $b(x,\phi_i)=\phi_ib(x)$ and the assumption that $\phi_i$ are Gaussian random variables
are crucial for availability of an explicit form of the likelihood
of the parameters of the random effects $\phi_i$. Indeed, assuming that 
$g(\varphi,\theta)d\nu(\varphi)\equiv N\left(\mu,\omega^2\right)$, \ctn{Maud12} obtain the likelihood
as the product of the following:
\begin{equation}
f_i(X_i|\theta)=\frac{1}{\left(1+\omega^2V_i\right)^{1/2}}\exp\left[-\frac{V_i}{2\left(1+\omega^2V_i\right)}
\left(\mu-\frac{U_i}{V_i}\right)^2\right]\exp\left(\frac{U^2_i}{2V_i}\right),
\label{eq:likelihood3}
\end{equation}
where $\theta=(\mu,\omega^2)\in\mathbb R\times\mathbb R^+$ ($R^+=(0,\infty)$),
and 
\begin{align}
&U_i=\int_0^{T_i}\frac{b(X_i(s))}{\sigma^2(X_i(s))}dX_i(s),\quad V_i=\int_0^{T_i}\frac{b^2(X_i(s))}{\sigma^2(X_i(s))}ds;
\quad i=1,\ldots,n,
\label{eq:sufficient}
\end{align}
are sufficient statistics. In (\ref{eq:likelihood3}), for $i=1,\ldots,n$, $X_i$ stands for 
$\left\{X_i(t);~t\in[0,T_i]\right\}$.

\ctn{Maud12} consider the $iid$ set-up by setting $x^i=x$ and $T_i=T$ for $i=1,\ldots,n$, and
directly prove weak consistency (convergence in probability) and
asymptotic normality of the $MLE$ of $\theta$. 
As an alternative, \ctn{Maitra14a} 
verify the regularity conditions of existing results in general set-ups provided in \ctn{Schervish95}
and \ctn{Hoadley71} to
prove asymptotic properties of the $MLE$ in this $SDE$ set-up. 
In the $iid$ set-up, this approach allowed \ctn{Maitra14a} to establish strong consistency of the $MLE$,
rather than weak consistency. Moreover, assumption (H4) of \ctn{Maud12}, requiring
that $b(\cdot)/\sigma(\cdot)$ is non-constant and for $i\geq 1$, $(U_i,V_i)$ admits a density
with respect to the Lebesgue measure on $\mathbb R\times \mathbb R^+$ which is jointly continuous 
and positive on an open ball of $\mathbb R\times\mathbb R^+$, was not required in their approach.
Also, not only in the $iid$ situation, \ctn{Maitra14a} 
prove asymptotic results related to the $MLE$ even in 
the independent but non-identical (we refer to this as non-$iid$) case.

To our knowledge, Bayesian asymptotics has not been investigated in the context of mixed effects models, even though
applied Bayesian analysis of such models is not rare (see, for example, \ctn{Wakefield94}, \ctn{Wakefield96}, \ctn{Bennett96}).
In this article, we consider the Bayesian framework associated with $SDE$-based random effects model, for both $iid$ and
non-$iid$ set-ups, and prove consistency and asymptotic normality of the Bayesian posterior distribution
of $\theta=(\mu,\omega^2)$. In other words, we consider prior distributions $\pi(\theta)$ of $\theta$ and study the
properties of the corresponding posterior
\begin{equation}
\pi_n(\theta|X_1,\ldots,X_n)=
\frac{\pi(\theta)\prod_{i=1}^nf_i(X_i|\theta)}
{\int_{\psi\in\Omega}\pi(\psi)\prod_{i=1}^nf_i(X_i|\psi)d\psi}
\label{eq:posterior1}
\end{equation}
as the sample size $n$ tends to infinity.
In what follows, in Section \ref{sec:consistency_iid} we investigate asymptotic properties of the 
posterior in the $iid$ context. In Section \ref{sec:consistency_non_iid} we investigate Bayesian
asymptotics in the non-$iid$ set-up.
We summarize our contribution and provide further discussion in Section \ref{sec:conclusion}.
Further details are provided in 
the supplement \ctn{Maitra14supp}, whose
sections, tables and figures have the prefix ``S-" when referred to in this paper.
Indeed in Section S-1 of the supplement we illustrate with examples when posterior consistency holds 
and fails; in the same section we also include examples pertaining to consistency and inconsistency
in a dependent set-up, which provide insights regarding extension of our asymptotic theory to 
dependent sets of $SDE$'s. 
In Section S-2 of the supplement we illustrate with an example the advantages of the Bayesian inference
over classical studies in $SDE$-based random effects models.
In Section S-3 we provide a brief discussion on choice of prior and associated posterior
computations in practical applications.

For the purpose of asymptotics, we adopt two further assumptions of \ctn{Maud12} (also adopted by \ctn{Maitra14a}), 
given by
\begin{itemize}
\item[(H3$^\prime$)] The parameter space $\Omega$ is a compact subset of $\mathbb R\times\mathbb R^+$.
\item[(H4$^\prime$)] The true value $\theta_0\in\Omega$.
\end{itemize}
Recall that condition (H4) of \ctn{Maud12}
was not required by \ctn{Maitra14a} in their classical approach. Neither is the
assumption required in our Bayesian approach. 
Also, as we shall show, our condition (H3$^\prime$) on compactness of $\Omega$, is not necessary for posterior consistency 
in the non-$iid$ case.
Notationally, ``$\stackrel{a.s.}{\rightarrow}$", ``$\stackrel{P}{\rightarrow}$" and ``$\stackrel{\mathcal L}{\rightarrow}$"
denote convergence ``almost surely", ``in probability" and ``in distribution", respectively.

\section{Consistency and asymptotic normality of the Bayesian posterior in the $iid$ set-up}
\label{sec:consistency_iid}

\subsection{Consistency of the Bayesian posterior distribution}
\label{subsec:Bayesian_consistency_iid}

Theorem 7.80 presented in \ctn{Schervish95} provides easy-to-verify sufficient conditions that 
ensure posterior consistency. We state the general theorem below, using which we prove posterior consistency
in our case.
\begin{theorem}[\ctn{Schervish95}]
\label{theorem:theorem3}
Let $\{X_n\}_{n=1}^{\infty}$ be conditionally $iid$ given $\theta$ with density $f_1(x|\theta)$
with respect to a measure $\nu$ on a space $\left(\mathcal X^1,\mathcal B^1\right)$. Fix $\theta_0\in\Omega$, and define,
for each $M\subseteq\Omega$ and $x\in\mathcal X^1$,
\[
Z(M,x)=\inf_{\psi\in M}\log\frac{f_1(x|\theta_0)}{f_1(x|\psi)}.
\]
Assume that for each $\theta\neq\theta_0$, there is an open set $N_{\theta}$ such that $\theta\in N_{\theta}$ and
that $E_{\theta_0}Z(N_{\theta},X_i)> -\infty$. 
Also assume that $f_1(x|\cdot)$ is continuous at $\theta$ 
for every $\theta$, a.s. $[P_{\theta_0}]$.
For $\epsilon>0$, define 
$C_{\epsilon}=\{\theta:\mathcal K_1(\theta_0,\theta)<\epsilon\}$, where 
\begin{equation}
\mathcal K_1(\theta_0,\theta)=E_{\theta_0}\left(\log\frac{f_1(X_1|\theta_0)}{f_1(X_1|\theta)}\right)
\label{eq:kl1}
\end{equation}
is the Kullback-Leibler divergence measure associated with observation $X_1$.
Let $\pi$ be a prior distribution such that $\pi(C_{\epsilon})>0$, for every $\epsilon>0$.
Then, for every $\epsilon>0$ and open set $\mathcal N_0$ containing $C_{\epsilon}$, the posterior satisfies
\begin{equation}
\lim_{n\rightarrow\infty}\pi_n\left(\mathcal N_0|X_1,\ldots,X_n\right)=1,\quad a.s.\quad [P_{\theta_0}].
\label{eq:posterior_consistency_iid}
\end{equation}
\end{theorem}

\subsubsection{Verification of posterior consistency}
\label{subsubsec:Bayesian_consistency_iid}
The condition $E_{\theta_0}Z(N_{\theta},X_i)> -\infty$ of the above theorem is 
verified in the context of Theorem 1 in 
\ctn{Maitra14a}. 
Here we provide a briefing on this. Note that in our case $f_1(x|\theta)$ 
is given by (\ref{eq:likelihood3}).
It then follows from the proof of Proposition 7 of \ctn{Maud12} that
for every $\theta\neq\theta_0$,
\begin{align}
\log\frac{f_1(x|\theta_0)}{f_1(x|\theta)}
&\geq -\frac{1}{2}\left\{\log\left(1+\frac{\omega^2}{\omega^2_0}\right)+\frac{|\omega^2-\omega^2_0|}{\omega^2}\right\}
-\frac{1}{2}|\omega^2_0-\omega^2|\left(\frac{U_1}{1+\omega^2_0V_1}\right)^2\left(1+\frac{\omega^2_0}{\omega^2}\right)\notag\\
&\quad -|\mu|\left\vert\frac{U_1}{1+\omega^2_0V_1}\right\vert\left(1+\frac{|\omega^2_0-\omega^2|}{\omega^2}\right)
-\left|\frac{\mu^2_0V_1}{2(1+\omega^2_0V_1)}\right|-\left|\frac{\mu_0U_1}{1+\omega^2_0V_1}\right|.\notag
\end{align}
Taking $N_{\theta}=\left(\underline\mu,\overline\mu\right)\times \left(\underline\omega^2,\overline\omega^2\right)$, and
making use of Lemma 1 of \ctn{Maud12} shows that $E_{\theta_0}\left(\frac{U_1}{1+\omega^2_0V_1}\right)^2$, 
$E_{\theta_0}\left\vert\frac{U_1}{1+\omega^2_0V_1}\right\vert$ and 
$E_{\theta_0}\left(\frac{\mu^2_0V_1}{2(1+\omega^2_0V_1)}\right)$
are finite. Hence, $E_{\theta_0}Z(N_{\theta},X_i)> -\infty$.

Now, all we need to ensure is that there exists a prior $\pi$ which gives positive probability to $C_{\epsilon}$
for every $\epsilon>0$. From the identifiability result given by Proposition 7 (i) of \ctn{Maud12}
it follows that $\mathcal K_1(\theta_0,\theta)=0$ if and only if $\theta=\theta_0$. Hence, for any $\epsilon>0$,
the set $C_{\epsilon}$ is non-empty, since it contains at least $\theta_0$.
In fact, since \ctn{Maud12} also show that
$\mathcal K_1(\theta_0,\theta)$ is continuous in $\theta$ (Proposition 7 (ii)), 
and since the parameter space $\Omega$ is compact, it follows that 
$\mathcal K_1(\theta_0,\theta)$ is uniformly continuous on $\Omega$. Hence,
for any $\epsilon>0$, there
exists $\delta_{\epsilon}$ which is independent of $\theta$, such that $\|\theta-\theta_0\|\leq\delta_{\epsilon}$ implies
$\mathcal K_1(\theta_0,\theta)<\epsilon$. In other words, 
$\left\{\theta:\|\theta-\theta_0\|\leq\delta_{\epsilon}\right\}\subseteq C_{\epsilon}$.

Let $\frac{d\pi}{d\nu}=h$ almost everywhere on $\Omega$, where $h(\theta)$ is any 
positive, continuous density on $\Omega$ with respect to the Lebesgue
measure $\nu$. By ``positive" density we mean a
density excluding any interval of null measure.
%
%
It then follows from the above arguments that
\begin{equation}
\pi(C_{\epsilon})\geq\pi\left(\left\{\theta:\|\theta-\theta_0\|\leq\delta_{\epsilon}\right\}\right)
\geq
\left[\underset{\left\{\theta:\|\theta-\theta_0\|\leq\delta_{\epsilon}\right\}}{\inf}~h(\theta)\right]\times
\nu\left(\left\{\theta:\|\theta-\theta_0\|\leq\delta_{\epsilon}\right\}\right)>0.
\label{eq:positivity1}
\end{equation}

Hence, (\ref{eq:posterior_consistency_iid}) holds in our case with any prior with positive, continuous density
with respect to the Lebesgue measure.
We summarize this result in the form of a theorem, stated below.
\begin{theorem}
\label{theorem:new_theorem3}
Assume the $iid$ set-up and conditions (H1$^\prime$) -- (H4$^\prime$). 
For $\epsilon>0$, define 
$C_{\epsilon}=\{\theta:\mathcal K_1(\theta_0,\theta)<\epsilon\}$, where 
\begin{align}
\mathcal K_1(\theta_0,\theta)
&=E_{\theta_0}\left(\log\frac{f_1(X_1|\theta_0)}{f_1(X_1|\theta)}\right)\notag\\
&=E_{\theta_0}\left[\frac{1}{2}\log\left(\frac{1+\omega^2V_1}{1+\omega^2_0V_1}\right)
+\frac{1}{2}\frac{(\omega^2_0-\omega^2)U^2_1}{(1+\omega^2V_1)(1+\omega^2_0V_1)}\right.\notag\\
&\left.\quad+\frac{\mu^2V_1}{2(1+\omega^2V_1)}-\frac{\mu U_1}{1+\omega^2V_1}
-\left(\frac{\mu^2_0V_1}{2(1+\omega^2_0V_1)}-\frac{\mu_0U_1}{1+\omega^2_0V_1}\right)\right]\notag\\
\label{eq:kl2}
\end{align}
is the Kullback-Leibler divergence measure associated with observation $X_1$.
Let the prior distribution $\pi$ of the parameter $\theta$ satisfy $\frac{d\pi}{d\nu}=h$ 
almost everywhere on $\Omega$, where $h(\theta)$ is any 
positive, continuous density on $\Omega$ with respect to the Lebesgue
measure $\nu$.  
Then the posterior (\ref{eq:posterior1}) is consistent
in the sense that for every $\epsilon>0$ and open set $\mathcal N_0$ containing $C_{\epsilon}$, the posterior satisfies
\begin{equation}
\lim_{n\rightarrow\infty}\pi_n\left(\mathcal N_0|X_1,\ldots,X_n\right)=1,\quad a.s.\quad [P_{\theta_0}].
\label{eq:posterior_consistency_iid2}
\end{equation}
\end{theorem}

\subsection{Asymptotic normality of the Bayesian posterior distribution}
\label{subsec:Bayesian_normality_iid}

We now investigate asymptotic normality of posterior distributions in our $SDE$ set-up. 
For our purpose, we make use of Theorem 7.102 in conjunction with 
Theorem 7.89 provided in \ctn{Schervish95}. 
These theorems make use of seven regularity conditions, of which only the first four
will be required for the $iid$ set-up. 
Hence, in this $iid$ context we state the four requisite conditions.

\subsubsection{Regularity conditions -- $iid$ case}
\label{subsubsec:regularity_iid}
\begin{itemize}
\item[(1)] The parameter space is $\Omega\subseteq\mathbb R^p$ for some finite $p$.
\item[(2)] $\theta_0$ is a point interior to $\Omega$.
\item[(3)] The prior distribution of $\theta$ has a density with respect to Lebesgue measure
that is positive and continuous at $\theta_0$.
\item[(4)] There exists a neighborhood $\mathcal N_0\subseteq\Omega$ of $\theta_0$ on which
$\ell_n(\theta)= \log f(X_1,\ldots,X_n|\theta)$ is twice continuously differentiable with 
respect to all co-ordinates of $\theta$, 
$a.s.$ $[P_{\theta_0}]$.
\end{itemize}

Before proceeding to justify asymptotic normality of our posterior, we furnish the relevant theorem below
(Theorem 7.102 of \ctn{Schervish95}).

\begin{theorem}[\ctn{Schervish95}]
\label{theorem:theorem4}
Let $\{X_n\}_{n=1}^{\infty}$ be conditionally $iid$ given $\theta$.
Assume the above four regularity conditions; 
also assume that there exists $H_r(x,\theta)$ such that, for each $\theta_0\in int(\Omega)$ and each
$k,j$,
\begin{align}
\sup_{\|\theta-\theta_0\|\leq r}\left\vert\frac{\partial^2}{\partial\theta_k\partial\theta_j}\log f_{1}(x|\theta_0)
-\frac{\partial^2}{\partial\theta_k\partial\theta_j}\log f_{1}(x|\theta)\right\vert\leq H_r(x,\theta_0),
\label{eq:H1}
\end{align}
with
\begin{equation}
\lim_{r\rightarrow 0}E_{\theta_0}H_r\left(X,\theta_0\right)=0.
\label{eq:H2}
\end{equation}
Further suppose that
the conditions of Theorem \ref{theorem:theorem3} hold, and that the Fisher's information matrix
$\mathcal I(\theta_0)$ is positive definite. 
Now denoting 
by $\hat\theta_n$ the $MLE$ associated with
$n$ observations, 
let
\begin{equation}
\Sigma^{-1}_n=\left\{\begin{array}{cc}
-\ell''_n(\hat\theta_n) & \mbox{if the inverse and}\ \ \hat\theta_n\ \ \mbox{exist}\\
\mathbb I_\tau & \mbox{if not},
\end{array}\right.
\label{eq:information1}
\end{equation}
where for any $t$,
\begin{equation}
\ell''_n(t)=\left(\left(\frac{\partial^2}{\partial\theta_i\partial\theta_j}\ell_n(\theta)\bigg\vert_{\theta=t}\right)\right),
\label{eq:information2}
\end{equation}
and $\mathbb I_\tau$ is the identity matrix of order $\tau$.
Thus, $\Sigma^{-1}_n$ is the observed Fisher's information matrix.

Letting $\Psi_n=\Sigma^{-1/2}_n\left(\theta-\hat\theta_n\right)$, for each compact subset 
$B$ of $\mathbb R^\tau$ and each $\epsilon>0$, the following holds:
\begin{equation}
\lim_{n\rightarrow\infty}P_{\theta_0}
\left(\sup_{\Psi_n\in B}\left\vert\pi_n(\Psi_n\vert X_1,\ldots,X_n)-\xi(\Psi_n)\right\vert>\epsilon\right)=0,
\label{eq:Bayesian_normality_iid}
\end{equation}
where $\xi(\cdot)$ denotes the density of the standard normal distribution.
\end{theorem}

\subsubsection{Verification of posterior normality}
\label{subsubsec:Bayesian_normality_iid}
We make the following assumption:
\begin{itemize}
\item[(H5$^\prime$)] The Fisher's information matrix $\mathcal I(\theta_0)$ is positive definite (see
\ctn{Maud12} for the form of $\mathcal I(\theta_0)$).
\end{itemize}
Now observe that the 
four regularity conditions in Section \ref{subsubsec:regularity_iid} trivially hold.
The remaining conditions of Theorem \ref{theorem:theorem4} are verified in the context of 
Theorem 2 of \ctn{Maitra14a}. Briefly, $\frac{\partial^2}{\partial\theta_k\partial\theta_j}\log f_{1}(x|\theta)$
is differentiable in $\theta=(\mu,\omega^2)$ and the derivative has finite expectation, which ensure (\ref{eq:H1})
and (\ref{eq:H2}).
Hence, (\ref{eq:Bayesian_normality_iid}) holds in our SDE set-up.
We summarize this result in the form of the following theorem.
\begin{theorem}
\label{theorem:new_theorem4}
Assume the $iid$ set-up and conditions (H1$^\prime$) -- (H5$^\prime$). Regarding (H3$^\prime$)
we assume, in particular, that $\theta_0\in int(\Omega)$.
Let the prior distribution $\pi$ of the parameter $\theta$ satisfy $\frac{d\pi}{d\nu}=h$ 
almost everywhere on $\Omega$, where $h(\theta)$ is any 
density with respect to the Lebesgue measure $\nu$ which is positive and continuous at $\theta_0$. 
Then, letting $\Psi_n=\Sigma^{-1/2}_n\left(\theta-\hat\theta_n\right)$, for each compact subset 
$B$ of $\mathbb R\times\mathbb R^+$ and each $\epsilon>0$, the following holds:
\begin{equation}
\lim_{n\rightarrow\infty}P_{\theta_0}
\left(\sup_{\Psi_n\in B}\left\vert\pi_n(\Psi_n\vert X_1,\ldots,X_n)-\xi(\Psi_n)\right\vert>\epsilon\right)=0.
\label{eq:Bayesian_normality_iid2}
\end{equation}
\end{theorem}

\section{Consistency and asymptotic normality of the Bayesian posterior in the non-$iid$ set-up}
\label{sec:consistency_non_iid}

In this section, as in \ctn{Maitra14a}, we do not enforce the restrictions $T_i=T$
and $x^i=x$ for $i=1,\ldots,n$. Consequently, here we deal with the set-up
where the processes $X_i(\cdot);~i=1,\ldots,n$, are independently,
but not identically distributed.
Following \ctn{Maitra14a}, we assume the following:
\begin{itemize}
\item[(H6$^\prime$)] The sequences $\{T_1,T_2,\ldots\}$ and 
$\{x^1,x^2,\ldots,\}$ are sequences in compact sets $\mathfrak T$ and $\mathfrak X$, respectively, so that
there exist convergent subsequences with limits in $\mathfrak T$ and $\mathfrak X$.
For notational convenience, we continue to denote the convergent subsequences as  $\{T_1,T_2,\ldots\}$
and $\{x^1,x^2,\ldots\}$. Let us denote the limits by $T^{\infty}$ and $x^{\infty}$, where $T^{\infty}\in\mathfrak T$
and $x^{\infty}\in\mathfrak X$. 
\end{itemize}


Following \ctn{Maitra14a}, we denote the process associated with the initial value $x$ and time point $t$ as $X(t,x)$, so that
$X(t,x^i)=X_i(t)$, and $X_i=\left\{X_i(t);~t\in[0,T_i]\right\}$. We also denote by
$\phi(x)$ the random effect
parameter associated with the initial value $x$ such that $\phi(x^i)=\phi_i$. 
We assume the following condition with respect to $\phi(x)$:
\begin{itemize}
\item[(H7$^\prime$)] $\phi(x)$ is a real-valued, continuous function of $x$, and that 
for $k\geq 1$, $\underset{x\in \mathfrak X}{\sup}~E\left[\phi(x)\right]^{2k}<\infty$. 
\end{itemize}
As in Proposition 1 of \ctn{Maud12}, assumption (H7$^\prime$) implies that for any $T>0$,
\begin{equation}
\underset{t\in [0,T],x\in \mathfrak X}{\sup}~E\left[X(t,x)\right]^{2k}<\infty.
\label{eq:sup_X}
\end{equation}

For $x\in \mathfrak X$ and $T\in \mathfrak T$, let
\begin{align}
U(x,T)&=\int_0^T\frac{b(X(s,x))}{\sigma^2(X(s,x))}d X(s,x)\label{eq:u_x_T};\\
V(x,T)&=\int_0^T\frac{b^2(X(s,x))}{\sigma^2(X(s,x))}ds.\label{eq:v_x_T}
\end{align}
For the non-$iid$ set-up we let $U_i=U(x^i,T_i)$ and $V_i=V(x^i,T_i)$.
As in \ctn{Maitra14a} we further assume that
\begin{itemize}
\item[(H8$^\prime$)]
\begin{equation}
\frac{b^2(x)}{\sigma^2(x)}<K(1+x^\tau),~\mbox{for some}~\tau\geq 1.
\label{eq:V_finite_moment2}
\end{equation}
\end{itemize}
This assumption ensures that moments of all orders of $V(x,T)$ are finite.
Then, by Theorem 3 of \ctn{Maitra14a}, the moments of uniformly integrable continuous 
functions of $U(x,T)$, $V(x,T)$ and $\theta$ are continuous in 
$x$, $T$ and $\theta$.
In particular, 
the Kullback-Leibler distance and the information matrix, which we denote by $\mathcal K_{x,T}(\theta_0,\theta)$
(or, $\mathcal K_{x,T}(\theta,\theta_0)$)
and $\mathcal I_{x,T}(\theta)$ respectively to emphasize
dependence on the initial values $x$ and $T$, are continuous in $x$, $T$ and $\theta$.
For $x=x^k$ and $T=T_k$, if we denote the Kullback-Leibler 
distance and the Fisher's information
as $\mathcal K_k(\theta_0,\theta)$ ($\mathcal K_k(\theta,\theta_0)$) and $\mathcal I_k(\theta)$, respectively,
then continuity of $\mathcal K_{x,T}(\theta_0,\theta)$ (or $\mathcal K_{x,T}(\theta,\theta_0)$) and $\mathcal I_{x,T}(\theta_0)$
with respect to $x$ and $T$ ensures that as $x^k\rightarrow x^{\infty}$ and $T_k\rightarrow T^{\infty}$,
$\mathcal K_{x^k,T_k}(\theta_0,\theta)\rightarrow \mathcal K_{x^{\infty},T^{\infty}}(\theta_0,\theta)=\mathcal K(\theta_0,\theta)$, say.
Similarly, $\mathcal K_{x^k,T_k}(\theta,\theta_0)\rightarrow \mathcal K(\theta,\theta_0)$ and 
$\mathcal I_{x^k,T_k}(\theta)\rightarrow \mathcal I_{x^{\infty},T^{\infty}}(\theta)=
\mathcal I(\theta)$, say. 
Thanks to compactness, the limits $\mathcal K(\theta_0,\theta)$, $\mathcal K(\theta,\theta_0)$ and $\mathcal I(\theta)$ are well-defined
Kullback-Leibler divergences and Fisher's information, respectively.
Consequently (see \ctn{Maitra14a}), the following hold for any $\theta\in\Omega$,
\begin{align}
\underset{n\rightarrow\infty}{\lim}~\frac{\sum_{k=1}^n\mathcal K_k(\theta_0,\theta)}{n}&=\mathcal K(\theta_0,\theta);
\label{eq:kl_limit_1}\\
\underset{n\rightarrow\infty}{\lim}~\frac{\sum_{k=1}^n\mathcal K_k(\theta,\theta_0)}{n}&=\mathcal K(\theta,\theta_0);
\label{eq:kl_limit_2}\\
\underset{n\rightarrow\infty}{\lim}~\frac{\sum_{k=1}^n\mathcal I_k(\theta)}{n}&=\mathcal I(\theta).
\label{eq:fisher_limit_1}
\end{align}
We assume that
\begin{itemize}
\item[(H9$^\prime$)] For any $\theta\in\Omega$, $\mathcal I(\theta)$ is positive definite.
\end{itemize}
The above results will be seen to have important roles as we proceed with the non-$iid$ Bayesian set-up.
For consistency in the Bayesian framework we utilize the theorem 
of \ctn{Choi04}, and for asymptotic normality of the posterior we make use of Theorem 7.89
of \ctn{Schervish95}.

\subsection{Posterior consistency in the non-$iid$ set-up}
\label{subsec:posterior_consistency_non_iid}

In our proceedings we need to ensure existence of moments
of the form
$$\underset{x\in\mathfrak X,T\in\mathfrak T}{\sup}~E_{\theta}\left[\exp\left\{\alpha \left|\omega^2_0-\omega^2\right| 
\left(\frac{U(x,T)}{1+\omega^2_0V(x,T)}\right)^2\left(1+\frac{\omega^2_0}{\omega^2}\right)\right\}\right],$$ for some
$0<\alpha<\infty$. The following extra assumption will be useful in this regard.
\begin{itemize}
\item[(H10$^\prime$)] There exists a strictly positive function $\alpha^*(x,T,\theta)$, 
continuous in $(x,T,\theta)$, such that for any $(x,T,\theta)$,
\begin{equation*}
E_{\theta}\left[\exp\left\{\alpha^*(x,T,\theta)K_1 U^2(x,T)\right\}\right]<\infty, 
\label{eq:moment_alpha1}
\end{equation*}
where
$K_1=\underset{\omega:~\theta\in\Omega}{\sup}~\left|\omega^2_0-\omega^2\right|\left(1+\frac{\omega^2_0}{\omega^2}\right)$.
\end{itemize}

Now, let 
\begin{equation}
\alpha^*_{\min}=\underset{x\in\mathfrak X,T\in\mathfrak T,\theta\in\Omega}{\inf}\alpha^*(x,T,\theta),
\label{eq:alpha_star}
\end{equation}
and 
\begin{equation}
\alpha=\min\left\{\alpha^*_{\min},c^*\right\},
\label{eq:alpha2}
\end{equation}
where $0<c^*<1/16$. 

Compactness ensures that $\alpha^*_{\min}>0$, so that $0<\alpha<1/16$.
It also holds due to compactness that for $\theta\in\Omega$,
\begin{equation}
\underset{x\in\mathfrak X,T\in\mathfrak T}{\sup}~E_{\theta}\left[\exp\left\{\alpha K_1U^2(x,T)\right\}\right]<\infty. 
\label{eq:moment_alpha2}
\end{equation}
This ensures that 
\begin{align}
&\underset{x\in\mathfrak X,T\in\mathfrak T}{\sup}~E_{\theta}\left[\exp\left\{\alpha \left|\omega^2_0-\omega^2\right| 
\left(\frac{U(x,T)}{1+\omega^2_0V(x,T)}\right)^2\left(1+\frac{\omega^2_0}{\omega^2}\right)\right\}\right]\notag\\
&\leq \underset{x\in\mathfrak X,T\in\mathfrak T}{\sup}~E_{\theta}\left[\exp\left\{\alpha K_1U^2(x,T)\right\}\right]\notag\\
&<\infty.
\label{eq:u_square_moment}
\end{align}

This choice of $\alpha$ ensuring (\ref{eq:moment_alpha2}) will be useful in verification
of the conditions of Theorem \ref{theorem:theorem5}, which we next state.

\begin{theorem}[\ctn{Choi04}]
\label{theorem:theorem5}
Let $\{X_i\}_{i=1}^{\infty}$ be independently distributed with densities $\{f_i(\cdot|\theta)\}_{i=1}^{\infty}$,
with respect to a common $\sigma$-finite measure, where $\theta\in\Omega$, a measurable space. The densities
$f_i(\cdot|\theta)$ are assumed to be jointly measurable. Let $\theta_0\in\Omega$ and let $P_{\theta_0}$
be the joint distribution of $\{X_i\}_{i=1}^{\infty}$ when $\theta_0$ is the true value of $\theta$.
Let $\{\Theta_n\}_{n=1}^{\infty}$ be a sequence of subsets of $\Omega$. Let $\theta$ have prior $\pi$ on $\Omega$.
Define the following:
\begin{align}
\Lambda_i(\theta_0,\theta) &=\log\frac{f_i(X_i|\theta_0)}{f_i(X_i|\theta)},\notag\\
\mathcal K_i(\theta_0,\theta) &= E_{\theta_0}\left(\Lambda_i(\theta_0,\theta)\right)\notag\\
\varrho_i(\theta_0,\theta) &= Var_{\theta_0}\left(\Lambda_i(\theta_0,\theta)\right).\notag
\end{align}
Make the following assumptions:
\begin{itemize}
\item[(1)] Suppose that there exists a set $B$ with $\pi(B)>0$ such that
\begin{enumerate}
\item[(i)] $\sum_{i=1}^{\infty}\frac{\varrho_i(\theta_0,\theta)}{i^2}<\infty,\quad\forall~\theta\in B$,
\item[(ii)] For all $\epsilon>0$, $\pi\left(B\cap\left\{\theta:\mathcal K_i(\theta_0,\theta)<\epsilon,~\forall~i\right\}\right)>0$.
\end{enumerate}
\item[(2)] Suppose that there exist test functions $\{\Phi_n\}_{n=1}^{\infty}$, sets $\{\Omega_n\}_{n=1}^{\infty}$
and constants $C_1,C_2,c_1,c_2>0$ such that
\begin{enumerate}
\item[(i)] $\sum_{n=1}^{\infty}E_{\theta_0}\Phi_n<\infty$,
\item[(ii)] $\underset{\theta\in \Theta^c_n\cap\Omega_n}{\sup}~E_{\theta}\left(1-\Phi_n\right)\leq C_1e^{-c_1n}$,
\item[(iii)] $\pi\left(\Omega^c_n\right)\leq C_2e^{-c_2n}$.
\end{enumerate}
\end{itemize}
Then,
\begin{equation}
\pi_n\left(\theta\in \Theta^c_n|X_1,\ldots,X_n\right)\rightarrow 0\quad a.s.~[P_{\theta_0}].
\label{eq:posterior_consistency_non_iid}
\end{equation}
\end{theorem}

\subsubsection{Validation of posterior consistency}
\label{subsubsec:posterior_consistency_non_iid}

Recall that $f_i(X_i|\theta)$ in our case is given by (\ref{eq:likelihood3}).
From the proof of Proposition 7 of \ctn{Maud12} it follows that 
$\left|\log \frac{f_i(X_i|\theta_0)}{f_i(X_i|\theta)}\right|$ has an upper bound which has finite expectation
and square of expectation under $\theta_0$, and is uniform for all 
$\theta\in B$, where $B$ is of the form $[\underline\mu,\overline\mu]\times[\underline\omega^2,\overline\omega^2]$, 
say, with $\underline\mu<\overline\mu$ and $0<\underline\omega^2<\overline\omega^2$. Hence, for each $i$, 
$\varrho_i(\theta_0,\theta)$ is finite. Moreover, since the sequences $\{T_1,T_2,\ldots\}$ and $\{x^1,x^2,\ldots\}$
belong to compact spaces $\mathfrak T$ and $\mathfrak X$, and the variance function $\varrho_{x,T}(\theta_0,\theta)$ viewed
as a function of $x$ and $T$, is bounded by a function
continuous in $x$ and $T$, $\varrho_i(\theta_0,\theta)<\kappa$, for some $0<\kappa<\infty$, uniformly in $i$.
Continuity of $\varrho_{x,T}(\theta_0,\theta)$ follows as an application of Theorem 3 of \ctn{Maitra14a}
where the required uniform integrability is assured by finiteness of the moments of all orders
of the random variable $U(x,T)/\left\{1+\omega^2V(x,T)\right\}$, 
for every $x\in\mathfrak X$, $T\in\mathfrak T$ (follows from Lemma 1 of \ctn{Maud12}) 
and compactness of $\mathfrak X$ and $\mathfrak T$.
Hence, choosing a prior that gives positive probability to the set $B$, it follows that for all $\theta\in B$,
$$ \sum_{i=1}^{\infty}\frac{\varrho_i(\theta_0,\theta)}{i^2}
<\kappa\sum_{i=1}^{\infty}\frac{1}{i^2}<\infty.$$ Hence, condition (1)(i) holds.

To verify (1)(ii) note that
%
because of compactness of $B$, $\mathcal K_i(\theta_0,\theta)$, which is continuous in $\theta$, is uniformly continuous
in $B$. Hence, for every $\epsilon>0$, there exists $\delta_i(\epsilon)$ independent of $\theta$ such that
$\|\theta-\theta_0\|<\delta_i(\epsilon)$ implies $\mathcal K_i(\theta_0,\theta)<\epsilon$.
Let us define 
\begin{equation}
\delta(\epsilon)=\inf\{\delta_{x,T}(\epsilon):~x\in \mathfrak X,T\in \mathfrak T\}, 
\label{eq:delta2}
\end{equation}
where $\delta_{x,T}(\epsilon)$ is any strictly positive continuous function of $x$ and $T$, depending upon $\epsilon$ such that
$\delta_{x^i,T_i}(\epsilon)=\delta_i(\epsilon)$, for every $i=1,2,\ldots$.
Compactness of $\mathfrak X$ and $\mathfrak T$ ensures that $\delta(\epsilon)>0$.
%
%
%
So, for any $\epsilon>0$, 
\begin{equation}
\left\{\theta\in B:\mathcal K_i(\theta_0,\theta)<\epsilon,~\forall~i\right\}
\supseteq\left\{\theta\in B:\|\theta-\theta_0\|<\delta(\epsilon)\right\}.
\label{eq:supset1}
\end{equation}
%
%
%
%
It follows that
\begin{equation}
\pi\left(B\cap\left\{\theta:\mathcal K_i(\theta_0,\theta)<\epsilon,~\forall~i\right\}\right)
\geq
\pi\left(\left\{\theta\in B:\|\theta-\theta_0\|<\delta(\epsilon)\right\}\right).
\label{eq:positivity2}
\end{equation}
The remaining part of the proof that the right hand side of (\ref{eq:positivity2}) is strictly positive,
follows exactly in the same way as the proof of strict positivity (\ref{eq:positivity1}) 
in Section \ref{subsubsec:Bayesian_consistency_iid}, with a positive, continuous prior density on $\Omega$
with respect to the Lebesgue measure.

We now verify conditions (2)(i), (2)(ii) and (2)(iii). We let $\Omega_n=\left(\Omega_{1n}\times\mathbb R^+\right)$,
where $\Omega_{1n}=\left\{\mu:|\mu|<M_n\right\}$, where $M_n=O(e^n)$. 
Note that
\begin{equation}
\pi\left(\Omega^c_n\right)=\pi\left(\Omega^c_{1n}\right)=\pi(|\mu|>M_n)
<E_{\pi}\left(|\mu|\right)M^{-1}_n,
\label{eq:sieve1}
\end{equation}
so that (2)(iii) holds, assuming that the prior $\pi$ is such that the expectation 
$E_{\pi}\left(|\mu|\right)$
is finite.

Fixing $\delta>0$, we construct the tests $\Phi_n$ as follows. 
\begin{align}
\Phi_n=\left\{\begin{array}{cc}
1 & \mbox{if}\quad \beta_n<\sqrt{e^{-n\delta}},\\
0 & \mbox{otherwise},\end{array}\right.
\end{align}
where 
\begin{equation}
\beta_n= \frac{L_n(\theta_0)}{\underset{\theta\in\Omega}{\sup}~L_n(\theta)}
=\frac{L_n(\theta_0)}{L_n(\hat\theta_n)}
\label{eq:beta_n}
\end{equation}
is the likelihood ratio test statistic under $H_0:\theta=\theta_0$ versus $H_1:\theta\neq\theta_0$.
Here $L_n(\theta)=\prod_{i=1}^nf_i(X_i|\theta)$, and, as before, $\hat\theta_n$ is the $MLE$ associated with
$n$ observations.
Now, denoting $-2\log\beta_n$ by $Z^2_n$, we obtain for $\alpha$ given by (\ref{eq:alpha2}),
\begin{eqnarray}
E_{\theta_0}\Phi_n&=P_{\theta_0}\left(\beta_n<\sqrt{e^{-n\delta}}\right)=P_{\theta_0}\left(Z^2_n>n\delta\right)
<e^{-\alpha n\delta}E_{\theta_0}\left(e^{\alpha Z^2_n}\right).
\label{eq:markov1}
\end{eqnarray}
Note that $Z^2_n=-n\left(\hat\theta_n-\theta_0\right)^T\frac{\ell''_n(\theta^*_n)}{n}\left(\hat\theta_n-\theta_0\right)$,
where $\ell_n(\theta)=\sum_{i=1}^n\log f_i(X_i|\theta)$, and $\theta^*_n$ lies between $\theta_0$ and
$\hat\theta_n$. 
Also, 
\begin{align}
\frac{\ell''_{n,ij}(\theta^*_n)}{n}&=\frac{\ell''_{n,ij}(\theta_0)}{n}
+(\theta^*_n-\theta_0)^T\frac{\ell'''_{n,ij}(\theta^{**}_n)}{n},
\label{eq:chisq0}
\end{align}
where $\ell''_{n,ij}$ is the $(i,j)$-th element of $\ell''_n$ and $\ell'''_{n,ij}$ is its derivative, and
$\theta^{**}_n$ lies between $\theta_0$ and $\theta^*_n$. Using Kolmogorov's strong law of large numbers for the non-$iid$ 
case (see, for example, \ctn{Serfling80}),
which holds in our problem due to finiteness of the moments of $U(x,T)/\left\{1+\omega^2V(x,T)\right\}$ 
for every $x$ and $T$ belonging to the compact
spaces $\mathfrak X$ and $\mathfrak T$, respectively, yields,
in conjunction with (\ref{eq:fisher_limit_1}), that  
\begin{equation}
\frac{\ell''_{n,ij}(\theta_0)}{n}\stackrel{a.s.}{\rightarrow}-\mathcal I_{ij}(\theta_0),
\label{eq:slln_ell_1}
\end{equation}
$\mathcal I_{ij}(\theta_0)$ being the $(i,j)$-th element of $\mathcal I(\theta_0)$.
Also, by Cauchy-Schwartz,
\begin{align}
\left|(\theta^*_n-\theta_0)^T\frac{\ell'''_{n,ij}(\theta^{**}_n)}{n}\right |
&\leq\|\theta^*_n-\theta_0\|\times\left\|\frac{\ell'''_{n,ij}(\theta^{**}_n)}{n}\right\|.
\label{eq:cs_ell_1}
\end{align}
In (\ref{eq:cs_ell_1}), due to boundedness of the third derivative (see the proof of Proposition 8
of \ctn{Maud12}), and due to continuity of the moments of $U(x,T)/\left\{1+\omega^2V(x,T)\right\}$ with respect to $x$ and $T$
(which follows from Theorem 3 of \ctn{Maitra14a} where uniform integrability is ensured by finiteness
of the moments of the aforementioned function for every $x,T$ belonging to compact sets $\mathfrak X$ and $\mathfrak T$), and then finally
applying Kolmogorov's strong law of large numbers for the non-$iid$ case, it can be easily shown that
$\left\|\frac{\ell'''_{n,ij}(\theta^{**}_n)}{n}\right\|=O_P(1)$. Since $\|\hat\theta_n-\theta_0\|=o_P(1)$, it follows that
$\|\theta^*_n-\theta_0\|=o_P(1)$ as well. Hence, 
\begin{equation}
\left|(\theta^*_n-\theta_0)^T\frac{\ell'''_{n,ij}(\theta^{**}_n)}{n}\right |=o_P(1),
\label{eq:ell_o_P_1}
\end{equation}
implying that
\begin{equation}
\frac{\ell''_n(\theta^*_n)}{n}\stackrel{P}{\rightarrow}-\mathcal I(\theta_0).
\label{eq:ell_2}
\end{equation}

Hence, due to (\ref{eq:ell_2}) and due to asymptotic normality 
of $MLE$ in our non-$iid$
set-up addressed in \ctn{Maitra14a}, 
under $P_{\theta_0}$, 
\begin{equation}
Z^2_n=-n\left(\hat\theta_n-\theta_0\right)^T\frac{\ell''_n(\theta^*_n)}{n}
\left(\hat\theta_n-\theta_0\right)
\stackrel{\mathcal L}{\rightarrow}\chi^2_1,
\label{eq:chisq}
\end{equation}
and so, by the continuous mapping theorem, $e^{\alpha Z^2_n}\stackrel{\mathcal L}{\rightarrow}e^{\alpha\chi^2_1}$.
Moreover, using the form $Z^2_n=-2\log\beta_n 
=-2\sum_{i=1}^n \left(\log f_i(X_i|\theta_0)-\log f_i(X_i|\hat\theta_n)\right)$,
we can write 
\begin{align}
E_{\theta_0}\left(e^{\alpha Z^2_n}\right)&=E_{X_1,\ldots,X_n|\theta_0}\left(e^{\alpha Z^2_n}\right)\\
&=E_{X_1,\ldots,X_n|\theta_0}\left[\exp\left\{-2\alpha\sum_{i=1}^n
\left(\log f_i(X_i|\theta_0)-\log f_i(X_i|\hat\theta_n)\right)\right\}\right]\notag\\
&=E_{\hat\theta_n|\theta_0}E_{X_1,\ldots,X_n|\theta_0,\hat\theta_n}\left[\exp\left\{-2\alpha\sum_{i=1}^n
\left(\log f_i(X_i|\theta_0)-\log f_i(X_i|\hat\theta_n)\right)\right\}\right]\notag\\
&=E_{\hat\theta_n|\theta_0}\prod_{i=1}^nE_{X_i|\theta_0,\hat\theta_n}\left[\exp\left\{-2\alpha
\left(\log f_i(X_i|\theta_0)-\log f_i(X_i|\hat\theta_n)\right)\right\}\right]\notag\\
&=E_{\hat\theta_n|\theta_0}\left[E_n(\hat\theta_n,\theta_0)\right],~\mbox{(say),}
\label{eq:ui_3}
\end{align}
where
\begin{align}
E_n(\hat\theta_n,\theta_0) &= \prod_{i=1}^nE_{X_i|\theta_0,\hat\theta_n}\left[\exp\left\{-2\alpha
\left(\log f_i(X_i|\theta_0)-\log f_i(X_i|\hat\theta_n)\right)\right\}\right].
\label{eq:E_n}
\end{align}

It follows from the lower bound obtained in the proof of Proposition 7 of \ctn{Maud12}, that 
conditional on $\hat\theta_n=\zeta=(\mu,\omega^2)$, 
$\log f_i(X_i|\theta_0)-\log f_i(X_i|\hat\theta_n)\geq C_3(U_i,V_i,\zeta)$,
where 
\begin{align}
C_3(U_i,V_i,\zeta)&= 
-\frac{1}{2}\left\{\log\left(1+\frac{\omega^2}{\omega^2_0}\right)+\frac{|\omega^2-\omega^2_0|}{\omega^2}\right\}
-\frac{1}{2}|\omega^2_0-\omega^2|\left(\frac{U_i}{1+\omega^2_0V_i}\right)^2\left(1+\frac{\omega^2_0}{\omega^2}\right)\notag\\
&\quad -|\mu|\left\vert\frac{U_i}{1+\omega^2_0V_i}\right\vert\left(1+\frac{|\omega^2_0-\omega^2|}{\omega^2}\right)
-\left|\frac{\mu^2_0V_i}{2(1+\omega^2_0V_i)}\right|-\left|\frac{\mu_0U_i}{1+\omega^2_0V_i}\right|.
\label{eq:lower_bound1}
\end{align}
Hence, for every given $n\geq 1$, due to the lower bound (\ref{eq:lower_bound1}) and 
assumption (H10$^\prime$), 
the latter implying 
(\ref{eq:u_square_moment}),
\begin{align}
E_n(\zeta,\theta_0) &= \prod_{i=1}^nE_{X_i|\theta_0,\hat\theta_n=\zeta}\left[\exp\left\{-2\alpha
\left(\log f_i(X_i|\theta_0)-\log f_i(X_i|\zeta)\right)\right\}\right]\notag\\
&\leq\prod_{i=1}^nE_{X_i|\theta_0,\hat\theta_n=\zeta}\left[\exp\left\{-2\alpha C_3(U_i,V_i,\zeta)\right\}\right]\notag\\
&<\infty,
\label{eq:E_n_finite}
\end{align}
for any $\zeta\in\Omega$.
Now, due to compactness of $\Omega$, 
$E_n(\zeta,\theta_0)\leq\underset{\vartheta\in\Omega}{\sup}~E_n(\vartheta,\theta_0)<\infty$, for every given $n$,
so that it follows from (\ref{eq:ui_3}), (\ref{eq:E_n}), (\ref{eq:lower_bound1}) and (\ref{eq:E_n_finite}) that 
$E_{\theta_0}\left(e^{\alpha Z^2_n}\right)= E_{\hat\theta_n|\theta_0}\left[E_n(\hat\theta_n,\theta_0)\right]
\leq \underset{\vartheta\in\Omega}{\sup}~E_n(\vartheta,\theta_0)<\infty$, for any given $n$.
So, for $n$ at most finite, 
\begin{equation}
\underset{n~\mbox{\tiny{at most finite}}}{\sup}~E_{\theta_0}\left(e^{\alpha Z^2_n}\right)<\infty.
\label{eq:finite1}
\end{equation}
In our problem, for large enough $n$, at most the following case can occur: 
for any given $\epsilon>0$, there exists $N_0(\epsilon)$ such that 
$\left|E_{\theta_0}\left(e^{\alpha Z^2_n}\right)
-E_{\theta_0}\left(e^{\alpha \chi^2_1}\right)\right|<\epsilon$ for $n\geq N_0(\epsilon)$,
where $E_{\theta_0}\left(e^{\alpha \chi^2_1}\right)<\infty$.
Combining this with (\ref{eq:finite1}) it follows that
\begin{equation}
\underset{n\geq 1}{\sup}~E_{\theta_0}\left(e^{\alpha Z^2_n}\right)<\infty.
\label{eq:finite2}
\end{equation}
Using this in conjunction with summation over (\ref{eq:markov1}), it is easily seen
that condition (2)(i) holds.

Let us now verify condition (2)(ii). For our purpose,
let us define $\Theta_n=\Theta_{\delta}=\left\{(\mu,\omega^2):\mathcal K(\theta,\theta_0)<\delta\right\}$,
where $\mathcal K(\theta,\theta_0)$, defined as in (\ref{eq:kl_limit_2}), is the proper Kullback-Leibler
divergence. 
Thus, $\mathcal K(\theta,\theta_0)>0$ if and only if $\theta\neq\theta_0$.
Now, 
\begin{align}
&E_{\theta}\left(1-\Phi_n\right)\notag\\
&=P_{\theta}\left(\beta_n>\sqrt{e^{-n\delta}}\right)
=P_{\theta}\left(-2\log\beta_n<n\delta\right)\notag\\
&=P_{\theta}\left(-n\left(\hat\theta_n-\theta\right)^T\frac{\ell''_n(\vartheta^*_n)}{n}\left(\hat\theta_n-\theta\right)
+2\ell_n(\theta_0)-2\ell_n(\theta)-2\left(\hat\theta_n-\theta\right)^T\ell'_n(\theta)>-n\delta\right)\notag\\
& \quad\quad\mbox{(here $\vartheta^*_n$ lies between $\theta$ and $\hat\theta_n$)}\notag\\
&<e^{\alpha n\delta}E_{\theta}\left(\exp\left\{-\alpha n\left(\hat\theta_n-\theta\right)^T\frac{\ell''_n(\vartheta^*_n)}{n}
\left(\hat\theta_n-\theta\right)
+2\alpha \ell_n(\theta_0)-2\alpha \ell_n(\theta)-2\alpha \left(\hat\theta_n-\theta\right)^T\ell'_n(\theta)\right\}\right),\notag\\
&\leq e^{\alpha n\delta}E_{\theta}\left(\exp\left\{-\alpha n\left(\hat\theta_n-\theta\right)^T\frac{\ell''_n(\vartheta^*_n)}{n}
\left(\hat\theta_n-\theta\right)
+2\alpha \ell_n(\theta_0)-2\alpha \ell_n(\theta)+2\alpha \left|\left(\hat\theta_n-\theta\right)^T\ell'_n(\theta)\right|\right\}\right)
\notag\\
&\leq e^{\alpha n\delta}\sqrt{E_{\theta}\left(\exp\left\{-2\alpha n\left(\hat\theta_n-\theta\right)^T\frac{\ell''_n(\vartheta^*_n)}{n}
\left(\hat\theta_n-\theta\right)
+4\alpha \ell_n(\theta_0)-4\alpha \ell_n(\theta)\right\}\right)}\label{eq:power0}\\
&\quad\quad\times \sqrt{E_{\theta}\left(\exp\left\{4\alpha \left|\left(\hat\theta_n-\theta\right)^T\ell'_n(\theta)\right|\right\}\right)}
\quad\quad\mbox{(using Cauchy-Schwartz inequality),}
\label{eq:power1}
\end{align}
where $\alpha$ is given by (\ref{eq:alpha2}).
Now observe that 
\begin{equation}
-n\left(\hat\theta_n-\theta\right)^T\frac{\ell''_n(\vartheta^*_n)}{n}
\left(\hat\theta_n-\theta\right)
\stackrel{\mathcal L}{\rightarrow}\chi^2_1,
\label{eq:chisq_1}
\end{equation}
\begin{equation}
\frac{\ell_n(\theta_0)- \ell_n(\theta)}{n}\stackrel{a.s.}{\rightarrow}-\mathcal K(\theta,\theta_0).
\label{eq:kl_1}
\end{equation}
The aforementioned convergence result (\ref{eq:kl_1}) is another
application of Kolmogorov's strong law of large numbers in the non-$iid$ case; as such,
similar arguments used to justify (\ref{eq:slln_ell_1}) remain valid here. 
%
%
The Cauchy-Schwartz inequality entails
\begin{align}
&\left|\left(\hat\theta_n-\theta\right)^T\ell'_n(\theta)\right|
=\left\vert(\hat\theta_n-\theta)^T{\mathcal I}(\theta)^{1/2}
{\mathcal I}(\theta)^{-1/2}\ell'_n(\theta)\right\vert\notag\\
&\leq \sqrt{n\left(\hat\theta_n-\theta\right)^T{\mathcal I}(\theta)\left(\hat\theta_n-\theta\right)}
\times \sqrt{n^{-1}\left\{\ell'_n(\theta)\right\}^T{\mathcal I}^{-1}(\theta)\ell'_n(\theta)},\notag\\
\label{eq:cs1}
\end{align}
where 
\begin{equation}
n\left(\hat\theta_n-\theta\right)^T{\mathcal I}(\theta)\left(\hat\theta_n-\theta\right)
\stackrel{\mathcal L}{\rightarrow}\chi^2_1
\label{eq:chisq2}
\end{equation}
and
\begin{equation}
n^{-1}\left\{\ell'_n(\theta)\right\}^T{\mathcal I}^{-1}(\theta)\ell'_n(\theta)
=n^{-1}tr~{\mathcal I}^{-1}(\theta)\ell'_n(\theta)\left\{\ell'_n(\theta)\right\}^T 
\stackrel{a.s.}{\rightarrow}tr\left({\mathcal I}^{-1}(\theta){\mathcal I}(\theta)\right)=2,
\label{eq:trace2}
\end{equation}
where, for any matrix $A$, $tr\left(A\right)$ denotes trace of the matrix $A$.

Hence, combining the asymptotic inequalities 
we obtain that for $\theta\in \Theta^c_n\cap\Omega_n$, where $n$ is sufficiently large,
\begin{align}
E_{\theta}\left(1-\Phi_n\right) &< 
e^{\alpha n\delta} 
\times e^{-2\alpha n\mathcal K(\theta,\theta_0)}
\times\sqrt{E_{\theta}\left(e^{2\alpha\chi^2_1}\right)\times E_{\theta}\left(e^{4\alpha\sqrt{2\chi^2_1}}\right)}\notag\\
&<
e^{\alpha n\delta}
\times e^{-2\alpha n\delta}\times
\sqrt{E_{\theta}\left(e^{2\alpha\chi^2_1}\right)\times E_{\theta}\left(e^{4\alpha\sqrt{2\chi^2_1}}\right)}\notag\\
&= e^{-\alpha\delta n}\times
\sqrt{E_{\theta}\left(e^{2\alpha\chi^2_1}\right)\times E_{\theta}\left(e^{4\alpha\sqrt{2\chi^2_1}}\right)}.
\label{eq:power2}
\end{align}
For our choice of $\alpha$, the expectations in (\ref{eq:power2}) are finite.
Also since the right hand side of (\ref{eq:power2}) does not depend upon $\theta$, (2)(ii) is proved in our case.
That is, finally, posterior consistency (\ref{eq:posterior_consistency_non_iid}) holds in our non-$iid$
SDE set-up. The result can be summarized in the form of the following theorem.
\begin{theorem}
\label{theorem:new_theorem5}
Assume the non-$iid$ SDE set-up. Also assume conditions (H1$^\prime$) and (H3$^\prime$) -- (H10$^\prime$).  
For any $\delta>0$, let $\Theta_{\delta}=\left\{(\mu,\omega^2):\mathcal K(\theta,\theta_0)<\delta\right\}$,
where $\mathcal K(\theta,\theta_0)$, defined as in (\ref{eq:kl_limit_2}), is the proper Kullback-Leibler
divergence. 
Let the prior distribution $\pi$ of the parameter $\theta$ satisfy $\frac{d\pi}{d\nu}=h$ 
almost everywhere on $\Omega$, where $h(\theta)$ is any 
positive, continuous density on $\Omega$ with respect to the Lebesgue
measure $\nu$.  
Then,
\begin{equation}
\pi_n\left(\theta\in \Theta^c_{\delta}|X_1,\ldots,X_n\right)\rightarrow 0\quad a.s.~[P_{\theta_0}].
\label{eq:posterior_consistency_non_iid2}
\end{equation}
\end{theorem}


\subsection{Asymptotic normality of the posterior distribution in the non-$iid$ set-up} 
\label{subsec:Bayesian_normality_non_iid}

For asymptotic normality of the posterior in the $iid$ situation, four regularity conditions,
stated in Section \ref{subsubsec:regularity_iid} were necessary. In the non-$iid$ framework,
three more are necessary, in addition to the already presented four conditions. They are as follows
(see \ctn{Schervish95} for details).

\subsubsection{Extra regularity conditions in the non-$iid$ set-up}
\label{subsubsec:regularity_non_iid}

\begin{itemize}
\item[(5)] The largest eigenvalue of $\Sigma_n$ goes to zero in probability.
\item[(6)] For $\delta>0$, define $\mathcal N_0(\delta)$ to be the open ball of radius $\delta$ around $\theta_0$.
Let $\rho_n$ be the smallest eigenvalue of $\Sigma_n$. If $\mathcal N_0(\delta)\subseteq\Omega$, there exists
$K(\delta)>0$ such that
\begin{equation}
\underset{n\rightarrow\infty}{\lim}~P_{\theta_0}\left(\underset{\theta\in\Omega\backslash\mathcal N_0(\delta)}{\sup}~
\rho_n\left[\ell_n(\theta)-\ell_n(\theta_0)\right]<-K(\delta)\right)=1.
\label{eq:extra1}
\end{equation}
\item[(7)] For each $\epsilon>0$, there exists $\delta(\epsilon)>0$ such that
\begin{equation}
\underset{n\rightarrow\infty}{\lim}~P_{\theta_0}\left(\underset{\theta\in\mathcal N_0(\delta(\epsilon)),\|\gamma\|=1}{\sup}~
\left\vert 1+\gamma^T\Sigma^{\frac{1}{2}}_n\ell''_n(\theta)\Sigma^{\frac{1}{2}}_n\gamma\right\vert<\epsilon\right)=1.
\label{eq:extra2}
\end{equation}
\end{itemize}

In the non-$iid$ case, the four regularity conditions presented in  Section \ref{subsubsec:regularity_iid}
and additional three provided above, are sufficient to guarantee (\ref{eq:Bayesian_normality_iid}).

\subsubsection{Verification of the regularity conditions}
\label{subsubsec:verify_last}

For $i=1,2$ and $j=1,2$, let the $(i,j)$-th element of $\ell''_n(\hat\theta_n)$ be denoted by
$\ell''_{n,ij}(\hat\theta_n)$. Then $\ell''_{n,ij}(\hat\theta_n)/n$ 
admits the following Taylor's series expansion around $\theta_0$:
\begin{equation}
\frac{\ell''_{n,ij}(\hat\theta_n)}{n}=\frac{\ell''_{n,ij}(\theta_0)}{n}+
\frac{(\hat\theta_n-\theta_0)^T\ell'''_{n,ij}(\theta^*_n)}{n},
\label{eq:last1}
\end{equation}
where $\theta^*_n$ lies between $\theta_0$ and $\hat\theta_n$.
In the same way as (\ref{eq:ell_2}) it can be easily shown that  
\begin{equation}
\frac{\ell''_n(\hat\theta_n)}{n}\stackrel{P}{\rightarrow}-\mathcal I(\theta_0).
\label{eq:ell_3}
\end{equation}
In other
words, $-\ell''_n(\hat\theta_n)$ and $n{\mathcal I}(\theta_0)$ are asymptotically equivalent (in probability).
Since the maximum eigenvalue of $n^{-1}{\mathcal I}^{-1}(\theta_0)$ goes to zero in probability
as $n\rightarrow\infty$,
so does the maximum eigenvalue of $\Sigma_n$.
Hence, condition (5) holds.

To verify condition (6), 
note that again by Kolmogorov's strong law of large numbers,
\begin{equation}
\frac{1}{n}\left(\ell_n(\theta)-\ell_n(\theta_0)\right)\stackrel{a.s.}{\rightarrow}-\mathcal K(\theta_0,\theta),
\label{eq:slln1}
\end{equation}
where, $\mathcal K(\theta_0,\theta)$ is given by (\ref{eq:kl_limit_1}),
%
%
Now, writing $\rho_n \left[\ell_n(\theta)-\ell_n(\theta_0)\right]$ as 
$n\rho_n \left[\frac{\ell_n(\theta)-\ell_n(\theta_0)}{n}\right]$ and
noting that $\Sigma_n=O_P\left(n^{-1}\right)$ implies $n\rho_n\stackrel{P}{\rightarrow}c$, where $c>0$, 
it follows from (\ref{eq:slln1}) that 
$\rho_n \left[\ell_n(\theta)-\ell_n(\theta_0)\right]\stackrel{P}{\rightarrow}-c\mathcal K(\theta_0,\theta)<0$.
Hence, condition (6) holds.

For condition (7) note that for $\theta\in\mathcal N_0(\delta(\epsilon))$, 
$\theta=\theta_0+\delta_2\frac{\theta_0}{\|\theta_0\|}$, where $0<\delta_2\leq\delta(\epsilon)$. 
So, using Taylor's series expansion around $\theta_0$, the $(i,j)$-th element of 
$\ell''_n(\theta)/n$ can be written as
\begin{equation}
\frac{\ell''_{n,ij}(\theta)}{n}=\frac{\ell''_{n,ij}(\theta_0)}{n}+
\delta_2\frac{\theta^T_0\ell'''_{n,ij}(\theta^*)}{n\|\theta_0\|},
\label{eq:last2}
\end{equation}
where $\theta^*$ lies between $\theta_0$ and $\theta$.
As $n\rightarrow\infty$, $\frac{\ell''_{n,ij}(\theta_0)}{n}$ tends, in probability,  
to the $(i,j)$-th element of $-{\mathcal I}(\theta_0)$. 
Now notice that
$$\frac{\left|\theta^T_0\ell'''_{n,ij}(\theta^*)\right|}{n\|\theta_0\|}
\leq \frac{\|\ell'''_{n,ij}(\theta^*)\|}{n},$$ 
so that $\frac{\left|\theta^T_0\ell'''_{n,ij}(\theta^*)\right|}{n\|\theta_0\|}=O_P(1)$ since
$\frac{\|\ell'''_{n,ij}(\theta^*)\|}{n}=O_P(1)$ as before.
Hence, it follows that $\ell''_n(\theta)=O_P\left(-n{\mathcal I}(\theta_0)+n\delta_2\right)$.
Since $\Sigma^{\frac{1}{2}}_n$ is asymptotically equivalent (in probability) to $n^{-\frac{1}{2}}
{\mathcal I}^{-\frac{1}{2}}(\theta_0)$, condition (7) holds.
We summarize our result in the form of the following theorem.
\begin{theorem}
\label{theorem:new_theorem6}
Assume the non-$iid$ set-up and conditions (H1$^\prime$) and (H3$^\prime$) -- (H9$^\prime$). Regarding (H3$^\prime$)
we assume, in particular, that $\theta_0\in int(\Omega)$.
Let the prior distribution $\pi$ of the parameter $\theta$ satisfy $\frac{d\pi}{d\nu}=h$ 
almost everywhere on $\Omega$, where $h(\theta)$ is any 
density with respect to the Lebesgue measure $\nu$ which is positive and continuous at $\theta_0$. 
Then, letting $\Psi_n=\Sigma^{-1/2}_n\left(\theta-\hat\theta_n\right)$, for each compact subset 
$B$ of $\mathbb R\times\mathbb R^+$ and each $\epsilon>0$, the following holds:
\begin{equation}
\lim_{n\rightarrow\infty}P_{\theta_0}
\left(\sup_{\Psi_n\in B}\left\vert\pi_n(\Psi_n\vert X_1,\ldots,X_n)-\xi(\Psi_n)\right\vert>\epsilon\right)=0.
\label{eq:Bayesian_normality_non_iid2}
\end{equation}
\end{theorem}

\section{Summary and discussion}
\label{sec:conclusion}


In this paper, we have investigated Bayesian posterior consistency 
in the context of $SDE$'s consisting of drift terms depending linearly upon 
random effect parameters. In particular, we have proved posterior consistency
and asymptotic normality in both $iid$ and non-$iid$ set-ups, as the number
of observed processes tends to infinity. 
In Section S-1 of our supplement we have illustrated our results with concrete examples, 
showing when posterior consistency will hold 
and not hold. Even in the dependent set-up we have illustrated, with examples, when consistency 
will hold and fail; see Section S-1.2 of the supplement. 
The latter examples can be looked upon as providing insights into the Bayesian 
asymptotic theory of dependent sets
of $SDE$'s.

It is also important to illustrate the value of Bayesian analysis in $SDE$-based random effects model,
particularly because, as per our results, at least asymptotically Bayesian analysis does not 
have edge over its classical counterpart. However, in small samples, Bayesian analysis can outperform
classical analysis when adequate prior knowledge on the parameter in question is available.
In Section S-2 of the supplement we present a simulation study to illustrate the advantage of 
Bayesian analysis in small samples. For realistic practical applications we include a brief 
discussion on elicitation of prior information and posterior computations in Section S-3 of the supplement.

Since discretization of the underlying continuous time processes is important for inference in 
practical situations, it is worth providing some remarks regarding discretization in our Bayesian
context. Firstly, note that discretized version of $U_i$ and $V_i$, as provided in \ctn{Maud12}, are
given by
\begin{align}
U^m_i&=\sum_{k=0}^{m-1}\frac{b(X_i(t_k))}{\sigma^2(X_i(t_k))}(X_i(t_{k+1})-X_i(t_k)),\notag\\
V^m_i&=\sum_{k=0}^{m-1}\frac{b^2(X_i(t_k))}{\sigma^2(X_i(t_k))}(t_{k+1}-t_k).\notag\\
\end{align}
Under mild conditions, Lemma 3 of \ctn{Maud12} provides bounds for the expectation associated with the differences
$U_i-U^m_i$ and $V_i-V^m_i$. We utilize the result to deduce that as $m\rightarrow\infty$,
$\pi_n(\cdot|X_1,\ldots,X_n)-\pi^{(m)}_n(\cdot|X_1,\ldots,X_n)=o_{P_{\theta_0}}(1)$, where 
$\pi^{(m)}_n(\cdot|X_1,\ldots,X_n)$ is the posterior density associated with $(U^m_i,V^m_i);~i=1,\ldots,n$.
Thus, it is not difficult to see that our consistency and asymptotic normality results for both $iid$ and non-$iid$ cases 
continue to hold as $m\rightarrow\infty$, $n\rightarrow\infty$ such that
$\frac{m}{n}\rightarrow\infty$.

Finally, it is important to remark that in this paper we have confined ourselves
to one-dimensional random effect parameters and one-dimensional $SDE$'s.
Although our non-$iid$ $SDE$ framework admits straightforward generalization to 
multi-dimensional situations (for the $iid$ counterpart multivariate generalization has been considered
by \ctn{Maud12}),
generalization of our asymptotic theory to high dimensions does not seem to be as straightforward. 
We reserve this problem for our future research.

\section*{Acknowledgments}
Sincere gratitude goes to three anonymous reviewers whose detailed comments have led to
much improved presentation of our article.
The first author gratefully acknowledges her CSIR Fellowship, Govt. of India.

\newpage

\renewcommand\thefigure{S-\arabic{figure}}
\renewcommand\thetable{S-\arabic{table}}
\renewcommand\thesection{S-\arabic{section}}

\setcounter{section}{0}
\setcounter{figure}{0}
\setcounter{table}{0}

\begin{center}
{\bf \Large Supplementary Material}
\end{center}

\section{Examples illustrating consistency and inconsistency}
\label{sec:examples}
In this section we consider simple examples in our $SDE$ set-up for illustrating
consistency and inconsistency in simple terms. For simplicity, in all the examples
we consider the set of $SDE$'s having the following form:
for $i=1,\ldots,n$,
\begin{equation}
d X_i(t)=\phi_idt+dW_i(t),
\label{eq:sde_example}
\end{equation}
that is, we set $b(\cdot)=\sigma(\cdot)\equiv 1$. Hence, for $i=1,\ldots,n$,
\begin{align}
U_i&=\phi_iT_i+W_i(T_i)~\mbox{and}~V_i=T_i.\label{eq:u_v_simple}
\end{align}


\subsection{Example 1}
\label{subsec:example1}
We assume that for $i=1,\ldots,n$, $\phi_i\stackrel{iid}{\sim}N\left(\mu,1\right)$, that is, we
set $\omega^2=1$. Letting $\mu_0$ be the true value of $\mu$, we investigate consistency 
of the posterior of $\mu$. Closed form expression of the posterior is available if we put $\pi(\mu)\equiv N\left(A,B^2\right)$
prior on $\mu$; here $-\infty<A<\infty$ and $B>0$. 
Indeed, $\pi_n(\mu|X_1,\ldots,X_n)\equiv N\left(\hat\mu_n,\hat\sigma^2_n\right)$, where
\begin{align}
\hat\mu_n &=\frac{\sum_{i=1}^n\frac{U_i}{1+V_i}+\frac{A}{B^2}}{\sum_{i=1}^n\frac{V_i}{1+V_i}+\frac{1}{B^2}};\label{eq:hat_mu}\\
\hat\sigma^2_n &=\frac{1}{\sum_{i=1}^n\frac{V_i}{1+V_i}+\frac{1}{B^2}}
=\frac{1}{\sum_{i=1}^n\frac{T_i}{1+T_i}+\frac{1}{B^2}}.\label{eq:hat_var}
\end{align}
Note that $\hat\mu_n$ follows the normal distribution with mean and variance given by the following: 
\begin{align}
E\left(\hat\mu_n\right)&=\frac{\mu_0\sum_{i=1}^n\frac{T_i}{1+T_i}+\frac{A}{B^2}}{\sum_{i=1}^n\frac{T_i}{1+T_i}+\frac{1}{B^2}};
\label{eq:mean_hat_mu}\\
Var\left(\hat\mu_n\right)&=\frac{\sum_{i=1}^n\frac{T_i}{1+T_i}}{\left(\sum_{i=1}^n\frac{T_i}{1+T_i}+\frac{1}{B^2}\right)^2}.
\end{align}
Now, for any $\epsilon>0$, by Chebychev's inequality,
\begin{align}
&P\left(|\hat\mu_n-\mu_0|>\epsilon\right)<\epsilon^{-4}E\left(\hat\mu_n-\mu_0\right)^4 
= O\left\{\left( Var\left(\hat\mu_n\right)\right)^2\right\}= 
O\left\{\left(\frac{\sum_{i=1}^n\frac{T_i}{1+T_i}}{\left(\sum_{i=1}^n\frac{T_i}{1+T_i}
+\frac{1}{B^2}\right)^2}\right)^2\right\}.
\label{eq:chebychev}
\end{align}

\subsubsection{Case 1: $T_i=T$ for all $i\geq 1$}
In this case, $Var(\hat\mu_n)=O\left(n^{-1}\right)$.
Hence, it follows from (\ref{eq:chebychev}) that 
for any $\epsilon>0$, $$\sum_{n=1}^{\infty}P\left(|\hat\mu_n-\mu_0|>\epsilon\right)<\infty,$$
so that $\hat\mu_n\stackrel{a.s.}{\rightarrow}\mu_0$. Moreover, since $\hat\sigma^2_n\rightarrow 0$, 
consistency of the posterior of $\mu$, for almost all data sequences, follows.

\subsubsection{Case 2: $T_i$ different and $\sum_{i=1}^{\infty}\frac{T_i}{1+T_i}=\infty$}
A typical instance of $T_i$ for which $\sum_{i=1}^{\infty}\frac{T_i}{1+T_i}=\infty$ is $\frac{T_i}{1+T_i}=\frac{1}{i}$
(so that $T_i=\frac{1}{i-1}$), for $i>1$; $T_1=c_0$, say, for any constant $c_0$.
%
For $\sum_{i=1}^{\infty}\frac{T_i}{1+T_i}=\infty$, $\hat\sigma^2_n\rightarrow 0$. Also, $E\left(\hat\mu_n\right)\rightarrow\mu_0$ and
$Var\left(\hat\mu_n\right)\rightarrow 0$, implying that $\hat\mu_n\stackrel{P}{\rightarrow}\mu_0$.
In fact, since $\sum_{i=1}^{n}\frac{T_i}{1+T_i}=O(n)$ as $n\rightarrow\infty$, it follows from (\ref{eq:chebychev}) that
for any $\epsilon>0$, $$\sum_{n=1}^{\infty}P\left(|\hat\mu_n-\mu_0|>\epsilon\right)<\infty.$$
Hence, consistency of the posterior of $\mu$ for almost all data sequences, follows.
Thus, in this situation, for the sequence $T_i=\frac{1}{i-1}$ for $i>1$; $T_1=c_0$ belonging to some
compact $\mathfrak T$ is appropriate.

Now observe that if $T_i\rightarrow\infty$ as $i\rightarrow\infty$, then also $\sum_{i=1}^{\infty}\frac{T_i}{1+T_i}=\infty$,
enforcing consistency.
Clearly in this case $\mathfrak T$ is non-compact. This example shows that compactness is sufficient, but not necessary
for consistency.

\subsubsection{Case 3: $\sum_{i=1}^{\infty}\frac{T_i}{1+T_i}<\infty$}
A typical instance of $T_i$ for which $\sum_{i=1}^{\infty}\frac{T_i}{1+T_i}<\infty$ is $\frac{T_i}{1+T_i}=\frac{1}{i^2}$
(so that $T_i=\frac{1}{i^2-1}$), for $i>1$; $T_1=c_0$, say, for any constant $c_0$.

In such cases, $\hat\sigma^2_n\nrightarrow 0$, showing that the posterior of $\mu$ is inconsistent.
In other words, even though the sequence $T_i=\frac{1}{i^2-1}$ for $i>1$; $T_1=c_0$ belongs to some
compact $\mathfrak T$, consistency still does not hold.
This example shows that not all convergent sequences in compact $\mathfrak T$ ensure consistency. As per our
theory, we can only assert that {\it there exists} at least one convergent subsequence of any sequence in $\mathfrak T$
for which consistency is attained.

\subsection{Example 2}
\label{subsec:example2}

In Example 1 we considered the set-up where, for $i=1,\ldots,n$, $\phi_i\stackrel{iid}{\sim}N\left(\mu,1\right)$.
This choice makes the $SDE$'s independent. Indeed, in this paper, we have provided theoretical results
assuming that the $SDE$'s are at least independent. However, at this point, we investigate with simulations,
when posterior consistency holds and fails, assuming a dependent set-up, the dependence induced by
the following $n$-variate normal distribution of $\bphi_n=(\phi_1,\ldots,\phi_n)$:
\begin{equation}
\bphi_n\sim N_n\left(\mu\bone_n,\omega^2\bSigma_n\right),
\label{eq:phi_mvn}
\end{equation}
where $\bone_n$ is the $n$-component vector with all entries 1, $\omega>0$, and $\bSigma_n$ is an $n\times n$ covariance
matrix. 
Let $\bU=(U_1,\ldots,U_n)'$,
and $\bV=diag\{V_1,\ldots,V_n\}$, that is, $\bV$ is an $n\times n$ diagonal matrix with the $i$-th diagonal element $V_i$.
Assuming that $\pi(\mu)\equiv N\left(A,B^2\right)$, it follows that the posterior of $\mu$ is given by
\begin{equation}
\pi_n(\mu|X_1,\ldots,X_n)\equiv N\left(\hat\mu_n,\hat\sigma^2_n\right),
\label{eq:dependent_posterior}
\end{equation}
where
\begin{align}
\hat\mu_n &=\frac{\omega^2\bone'\bSigma^{-1}_n\left(\omega^2\bV+\bSigma^{-1}_n\right)\bU+\frac{A}{B^2}}
{\bone'\bSigma^{-1}_n\left\{\bSigma_n-\left(\omega^2\bV+\bSigma^{-1}_n\right)^{-1}\right\}\bSigma^{-1}_n\bone
+\frac{1}{B^2}};\label{eq:dependent_hat_mu}\\
\hat\sigma^2_n &=\frac{\omega^2}{\bone'\bSigma^{-1}_n\left\{\bSigma_n-\left(\omega^2\bV+\bSigma^{-1}_n\right)^{-1}\right\}
\bSigma^{-1}_n\bone+\frac{1}{B^2}}.\label{eq:dependent_hat_var}
\end{align}
We conduct two simulation studies to investigate consistency in this set-up.

\subsubsection{First simulation study -- weak dependence structure}
\label{subsubsec:first_simstudy}
As in Example 1, we set $b(\cdot)=\sigma(\cdot)\equiv 1$ in our $SDE$'s, so that 
$U_i=\phi_iT_i+W_i(T_i)$ and $V_i=T_i$; we set $T_i=5$ for each $i$. For the distribution of $\bphi_n$ 
of the form (\ref{eq:phi_mvn}) we set the true value $\mu_0=1$. Also, for simplicity, we set $\omega^2=1$. 
For the covariance matrix $\bSigma_n$, we consider the following weakly dependent structure: 
\begin{equation}
\bSigma_n=\left(\begin{array}{ccccccc}
1 & \frac{1}{3} & 0 & 0 & 0 &\cdots & 0\\
\frac{1}{3} & 1 & \frac{1}{3} & 0 & 0 & \cdots & 0\\
0 & \frac{1}{3} & 1 & \frac{1}{3} & 0 & \cdots & 0\\
\vdots & \vdots & \vdots & \vdots & \vdots & \vdots & \vdots\\
0 & 0 & 0 & \cdots & 0 & \frac{1}{3} & 1\\
\end{array}
\right).
\label{eq:weak_dependence}
\end{equation}
For any $n$, this matrix is strictly diagonally dominant, and hence positive definite. 

We generate the data by generating $\bphi_n$ from the $n$-variate normal
(\ref{eq:phi_mvn}), generating $W_i(T_i)$ and then forming $U_i$ and $V_i$, for $i=1,\ldots,n$. 

\begin{figure}
\centering
\includegraphics[height=7cm,width=7cm]{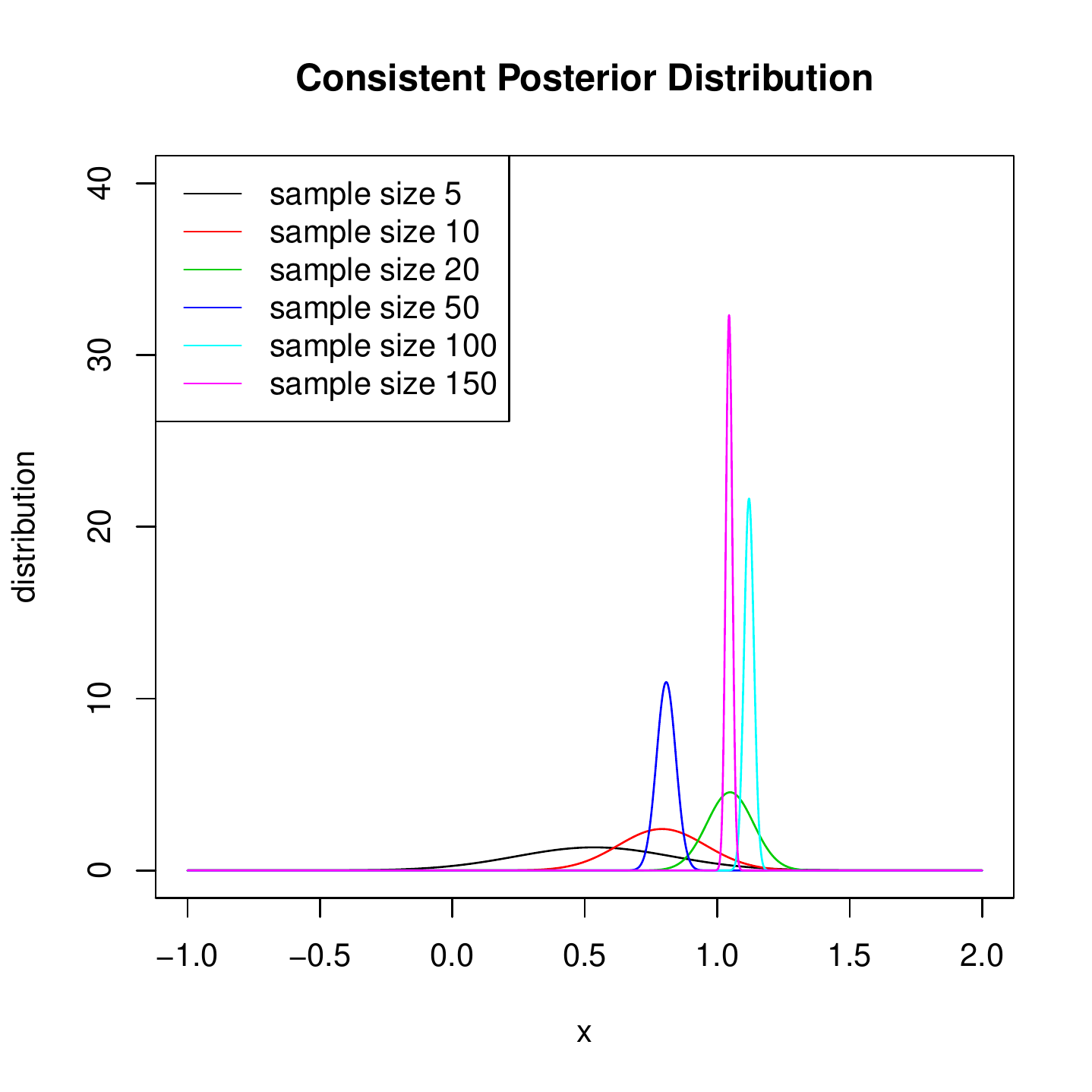}
\caption{Illustration of posterior consistency in the weakly dependent set-up.}
\label{fig:consistent}
\end{figure}

Figure \ref{fig:consistent} displays the posterior distribution of $\mu$ for various
sample sizes. Clearly, for larger sample sizes, the posterior increasingly concentrates 
around the true value $\mu_0=1$, thus demonstrating posterior consistency.

\subsubsection{Second simulation study -- strong dependence structure}
\label{subsubsec:second_simstudy}
In the second experiment, we consider a strong dependence structure between the components of $\bphi_n$,
quantified by 
\begin{equation}
\bSigma_n=\left(\begin{array}{ccccccc}
1 & \frac{1}{3} & \frac{1}{3} & \frac{1}{3} & \frac{1}{3} &\cdots & \frac{1}{3}\\
\frac{1}{3} & 1 & \frac{1}{3} & \frac{1}{3} & \frac{1}{3} & \cdots & \frac{1}{3}\\
\frac{1}{3} & \frac{1}{3} & 1 & \frac{1}{3} & \frac{1}{3} & \cdots & \frac{1}{3}\\
\vdots & \vdots & \vdots & \vdots & \vdots & \vdots & \vdots\\
\frac{1}{3} & \frac{1}{3} & \frac{1}{3} & \cdots & \frac{1}{3} & \frac{1}{3} & 1\\
\end{array}
\right).
\label{eq:strong_dependence}
\end{equation}
The rest of the set-up remains exactly same as in the previous experiment with the weak dependence structure.

\begin{figure}
\centering
\includegraphics[height=7cm,width=7cm]{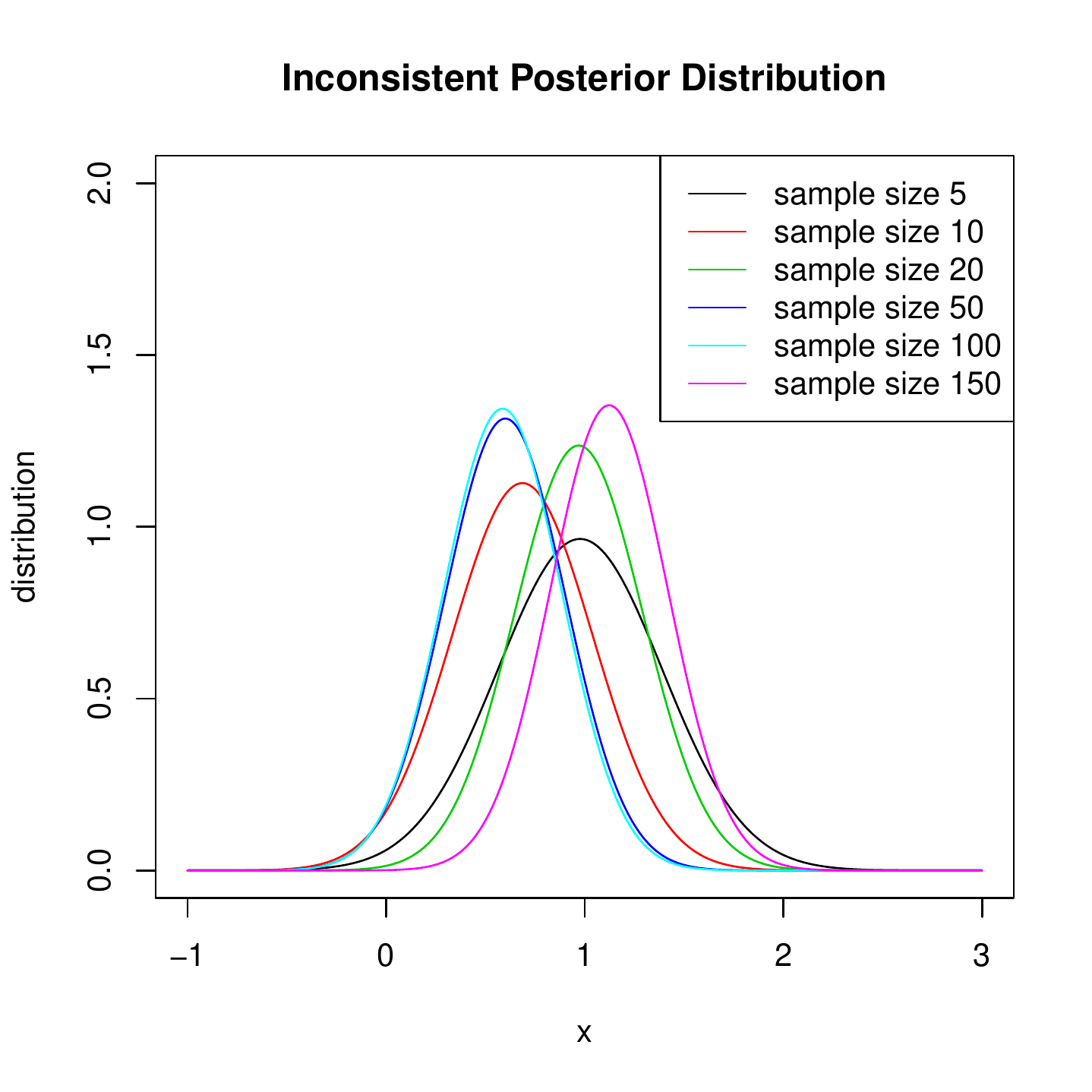}
\caption{Illustration of posterior inconsistency in the strongly dependent set-up.}
\label{fig:inconsistent}
\end{figure}

Figure \ref{fig:inconsistent} displays the posterior distribution of $\mu$ for various
sample sizes, in the strongly dependent situation. Note that even for large sample sizes, the 
variability of the posterior distribution does not decrease, which is clearly indicative
of posterior inconsistency.
Hence, in the strongly dependent $SDE$ set-up, posterior consistency does not hold. 

\section{Illustration of advantages of Bayesian analysis over classical inference in small samples}
\label{sec:bayes_advantages}

It is important to illustrate the value of Bayesian analysis in $SDE$-based random effects model,
particularly because, as per our results, at least asymptotically Bayesian analysis does not 
have edge over its classical counterpart. However, in small samples, Bayesian analysis can outperform
classical analysis when adequate prior knowledge on the parameter in question is available.
We undertake a simulation study to illustrate the advantage of Bayesian analysis in small samples. Briefly,
we consider the same set-up as Section \ref{sec:examples}, assuming that the true value $\mu_0=1$. 
We choose the prior $\pi(\mu)\equiv N\left(A,B^2\right)$ with $A=0$ and $B=1.5$. 
In this set-up, we obtain the 95\% confidence interval in the classical case and 95\% highest
posterior density (HPD) interval in the Bayesian case, for the true parameter $\mu_0$.
The results are presented in Table \ref{table:comparison}. Note that all the intervals, for both
classical and Bayesian analyses include the true value $\mu_0=1$, however, the lengths of the Bayesian
95\% HPD intervals are significantly shorter than the corresponding 95\% classical confidence intervals
in all the cases. In fact, lesser the sample size, larger is the difference between the lengths
of the Bayesian and classical intervals. In real situations involving random effects, adequately informative 
prior opinions regarding the parameter in question can often be obtained, and this small example demonstrates
that such information can substantially enhance Bayesian inference. 
\begin{table}[h]
\centering
\caption{Comparison between classical and Bayesian analysis of $SDE$-based random effects
model for small samples.}
\label{table:comparison}
\begin{tabular}{|c||c|c||c|c|}
\hline
Sample & Classical 95\% Confidence & Interval & Bayes 95\% HPD & Interval\\
Size & Interval & Length & Interval & Length\\
\hline
2 & (-1.052, 3.243) & 4.294 & (-0.934, 1.866) & 2.801\\ 
3 & (-0.744, 2.293) & 3.036 & (-0.716, 1.681) & 2.397\\
4 & (-0.607, 1.872) & 2.479 & (-0.597, 1.508) & 2.105\\
5 & (-0.526, 1.621) & 2.147 & (-0.521, 1.374) & 1.894\\
6 & (-0.470, 1.450) & 1.920 & (-0.468, 1.268) & 1.735\\
7 & (-0.429, 1.324) & 1.753 & (-0.428, 1.182) & 1.610\\
8 & (-0.397, 1.226) & 1.623 & (-0.397, 1.112) & 1.508\\
9 & (-0.372, 1.146) & 1.518 & (-0.371, 1.052) & 1.423\\
10 & (-0.351, 1.081) & 1.431 & (-0.350, 1.001) & 1.351\\
\hline
\end{tabular}
\end{table}

\section{Brief discussion on prior elicitation and posterior computation in realistic situations}
\label{sec:prior_computation}

In our simulation studies, for simplicity of illustrations, we have considered the $N\left(A,B^2\right)$ prior
on $\mu$, assuming $A$ and $B$ to be known; also we have set $\omega^2=1$. In realistic situations elicitation
of such strong prior information is not always straightforward. 
Moreover, the random effects parameters $\phi_i$ may be $d$-dimensional, 
so that $\mu$ is to be replaced by the $d$-dimensional vector $\bmu=\left(\mu_1,\ldots,\mu_d\right)$,
and $\omega^2$ needs to be replaced with the $d\times d$ matrix $\bSigma$.
The multidimensional situation makes appropriate choices of priors
even more difficult. In the context of $SDE$-based pharmacokinetic models \ctn{Yan14}, following \ctn{Cruz06},
proposed independent normal priors for the components of $\bmu$, given by
$\mu_i\sim N\left(A_i,B^2_i\right)$, 
and an inverse Wishart prior for $\bSigma$ with scale matrix $\bSigma_0$ and degrees of freedom $d+1$.
However, because of the difficulties of eliciting information, they assume
non-informative priors for the hyperparameters.

Informative priors can be elicited if historical data, that is, data associated with previous studies, are available. 
Then, following \ctn{Yan14}, using the aforementioned non-informative priors, one can first obtain
the posterior distributions of $\bmu$ and $\bSigma$, given only the historical data. These posterior distributions
based on historical data can then be used as informative prior distributions for Bayesian analysis of the current data.
Such an approach has been advocated and formalized in the context of generalized linear mixed models
by \ctn{Ibrahim00}.

It is important to remark that for complicated priors as discussed above, the posterior distribution
need not be available in closed form. Even obtaining closed form of the prior, which is the posterior given the
historical data, is not guaranteed. In the context of generalized linear mixed effects model,
\ctn{Ibrahim00} propose a Gibbs sampling algorithm based on the centering
strategies of \ctn{Gelfand96a}, for efficiently sampling from the desired posterior, given the current data.
However, since our $SDE$-based model does not fall within the class of generalized linear mixed models, 
such Gibbs sampling strategies need not be available in our case. Instead, the Metropolis-Hastings (MH) method
can be used to sample from the posterior. In the context of pharmacokinetic models, \ctn{Yan14}
estimate the parameters using a combination of extended Kalman filter and random walk MH algorithm;
see also \ctn{Donnet13} for a comprehensive review on various techniques of classical and Bayesian 
estimation of $SDE$'s for pharmacokinetic/pharmacodynamic
models.

For a fully Bayesian approach, required posterior computations can be carried out using only the MH methodology.
The MH methodology, however, has drawbacks in that convergence can often be quite slow and computations
can be very burdensome, when the dimension $d$ is large. To bypass these problems, \ctn{Dutta14}
have developed a novel methodology which they refer to as Transformation based Markov Chain Monte Carlo 
(TMCMC), that can update the entire
high-dimensional parameter set in a single block using simple deterministic transformations of a single
random variable, thus effectively reducing the multidimensional parameter to a single dimension. Apart from
drastically reducing computing time, the method promises much improved acceptance rates and convergence properties; see
\ctn{Dutta14}, \ctn{Dey15a}, \ctn{Dey15b}, \ctn{Dey15c}.

\normalsize
\bibliographystyle{natbib}
\bibliography{irmcmc}

\end{document}